\documentclass[11pt]{article}
\usepackage{amsmath,amssymb}
\usepackage{graphicx, citesort}

\def\sgn{{\hbox{sgn}}}

\def\R{{{\mathbf{R}}}} 

\def\P{{\hbox{\bf P}}}

\newenvironment{proof}{\noindent {\bf Proof} }{\endprf\par}
\def \endprf{\hfill {\vrule height6pt width6pt depth0pt}\medskip}
\def\emph#1{{\it #1}} \def\textbf#1{{\bf #1}}

\newcommand{\<}{\langle}
\renewcommand{\>}{\rangle}

\newcommand{\goto}{\rightarrow}

\newtheorem{theorem}{Theorem}[section]
\newtheorem{lemma}[theorem]{Lemma}
\newtheorem{corollary}[theorem]{Corollary}

\newtheorem{definition}[theorem]{Definition}

\numberwithin{equation}{section}

\def\cF{{\cal F}}

\begin{document}

\title{Decoding by Linear Programming}

\author{Emmanuel Candes$^{\dagger}$ and Terence Tao$^{\sharp}$\\
  \vspace{-.3cm}\\
  $\dagger$ 
Applied and Computational Mathematics, Caltech, Pasadena, CA 91125\\
  \vspace{-.3cm}\\
  $\sharp$ 
Department of Mathematics, University of California, Los Angeles, CA 90095
} 

\date{December 2004} 

\maketitle

\begin{abstract}

  This paper considers the classical error correcting problem which is
  frequently discussed in coding theory. We wish to recover an input
  vector $f \in \R^n$ from corrupted measurements $y = A f + e$. Here,
  $A$ is an $m$ by $n$ (coding) matrix and $e$ is an arbitrary and
  unknown vector of errors. Is it possible to recover $f$ exactly from
  the data $y$?
  
  We prove that under suitable conditions on the coding matrix $A$,
  the input $f$ is the unique solution to the $\ell_1$-minimization
  problem ($\|x\|_{\ell_1} := \sum_i |x_i|$)
  $$  
  \min_{g \in \R^n} \,\, \| y - Ag \|_{\ell_1} 
  $$
  provided that the support of the vector of errors is not too
  large, $\|e\|_{\ell_0} := |\{i : e_i \neq 0\}| \le \rho \cdot m$ for
  some $\rho > 0$. In short, $f$ can be recovered exactly by solving a
  simple convex optimization problem (which one can recast as a linear
  program).  In addition, numerical experiments suggest that this
  recovery procedure works unreasonably well; $f$ is recovered exactly
  even in situations where a significant fraction of the output is corrupted.
  
  This work is related to the problem of finding sparse solutions to
  vastly underdetermined systems of linear equations. There are also
  significant connections with the problem of recovering signals from
  highly incomplete measurements. In fact, the results introduced in
  this paper improve on our earlier work \cite{OptimalRecovery}.
  Finally, underlying the success of $\ell_1$ is a crucial property we
  call the uniform uncertainty principle that we shall describe in
  detail.

\end{abstract}


{\bf Acknowledgments.} E.~C. is partially supported by National
Science Foundation grants DMS 01-40698 (FRG) and ACI-0204932 (ITR),
and by an Alfred P.  Sloan Fellowship. T.~T. is supported in part by a
grant from the Packard Foundation.  Many thanks to Rafail Ostrovsky
for pointing out possible connections between our earlier work and the
decoding problem.  E.~C.~would also like to acknowledge inspiring
conversations with Leonard Schulmann and Justin Romberg.

\pagebreak

\section{Introduction}
\label{sec:introduction}

\subsection{Decoding of linear codes}

This paper considers the model problem of recovering an input vector
$f \in \R^n$ from corrupted measurements $y = A f + e$. Here, $A$ is
an $m$ by $n$ matrix (we will assume throughout the paper that $m >
n$), and $e$ is an arbitrary and unknown vector of errors. The problem
we consider is whether it is possible to recover $f$ exactly from the
data $y$. And if so, how? 

In its abstract form, our problem is of course equivalent to the
classical error correcting problem which arises in coding theory as we
may think of $A$ as a \emph{linear code}; a linear code is a given
collection of codewords which are vectors $a_1,\ldots, a_n \in
\R^m$---the columns of the matrix $A$. Given a vector $f \in \R^n$ (the
``plaintext'') we can then generate a vector $Af$ in $\R^m$ (the
``ciphertext''); if $A$ has full rank, then one can clearly recover
the plaintext $f$ from the ciphertext $Af$. But now we suppose that
the ciphertext $Af$ is corrupted by an arbitrary vector $e \in \R^m$
giving rise to the corrupted ciphertext $Af+e$.  The question is then:
given the coding matrix $A$ and $Af+e$, can one recover $f$ exactly?

As is well-known, if the fraction of the corrupted entries is too
large, then of course we have no hope of reconstructing $f$ from
$Af+e$; for instance, assume that $m = 2n$ and consider two distinct
plaintexts $f, f'$ and form a vector $g \in \R^m$ by setting half of
its $m$ coefficients equal to those of $Af$ and half of those equal to
those of $A f'$.  Then $g = Af + e = Af' + e'$ where both $e$ and $e'$
are supported on sets of size at most $n = m/2$. This simple example
shows that accurate decoding is impossible when the size of the
support of the error vector is greater or equal to a half of that of
the output $A f$. Therefore, a common assumption in the literature is
to assume that only a small fraction of the entries are actually
damaged
\begin{equation}
  \label{eq:l0}
\|e\|_{\ell_0} := |\{i : e_i \neq 0\}| \le \rho \cdot m.  
\end{equation}
For which values of $\rho$ can we hope to reconstruct $e$ with
practical algorithms? That is, with algorithms whose complexity is at
most polynomial in the length $m$ of the code $A$?

To reconstruct $f$, note that it is obviously sufficient to
reconstruct the vector $e$ since knowledge of $Af+e$ together with $e$
gives $Af$, and consequently $f$ since $A$ has full rank.  Our
approach is then as follows. We construct a matrix which annihilates
the $m \times n$ matrix $A$ on the left, i.e. such that $F A = 0$.
This can be done in an obvious fashion by taking a matrix $F$ whose
kernel is the range of $A$ in $\R^m$, which is an $n$-dimensional
subspace (e.g. $F$ could be the orthogonal projection onto the
cokernel of $A$). We then apply $F$ to the output $y = A f + e$ and
obtain
\begin{equation}
  \label{eq:Be}
  \tilde y = F(Af + e) = F e
\end{equation}
since $F A = 0$.  Therefore, the decoding problem is reduced to that of
reconstructing a \emph{sparse} vector $e$ from the observations $F e$
(by sparse, we mean that only a fraction of the entries of $e$ are
nonzero).

\subsection{Sparse solutions to underdetermined systems}

Finding sparse solutions to underdetermined systems of linear
equations is in general $NP$-hard \cite{Natarajan,Equivl0l1}. For
example, the sparsest solution is given by 
\begin{equation}
  \label{eq:P0}
  (P_0) \quad \quad 
\min_{d \in \R^m} \|d\|_{\ell_0} \quad \text{subject to} 
\quad F d = \tilde y \,\, ( = F e), 
\end{equation}
and to the best of our knowledge, solving this problem essentially
require exhaustive searches over all subsets of columns of $F$, a
procedure which clearly is combinatorial in nature and has exponential
complexity.

This computational intractability has recently led researchers to
develop alternatives to $(P_0)$, and a frequently discussed approach
considers a similar program in the $\ell_1$-norm which goes by
the name of {\em Basis Pursuit} \cite{BP}:
\begin{equation}
\label{eq:(P_1)}
(P_1)\quad \quad  \min_{x \in \R^m} 
\|d\|_{\ell_1}, \qquad F d = \tilde y, 
\end{equation}
where we recall that $\|d\|_{\ell_1} = \sum_{i = 1}^m |d_i|$.  Unlike
the $\ell_0$-norm which enumerates the nonzero coordinates, the
$\ell_1$-norm is convex. It is also well-known \cite{Bloomfield} that
$(P_1)$ can be recast as a linear program (LP).

Motivated by the problem of finding sparse decompositions of special
signals in the field of mathematical signal processing and following
upon the ground breaking work of Donoho and Huo \cite{DonohoHuo}, a
series of beautiful articles
\cite{GribonvalNielsen,DonohoElad,EladBruckstein,Tropp03} showed exact
equivalence between the two programs $(P_0)$ and $(P_1)$. In a
nutshell, this work shows that for $m/2$ by $m$ matrices $F$ obtained
by concatenation of two orthonormal bases, the solution to both
$(P_0)$ and $(P_1)$ are unique and identical provided that in the most
favorable case, the vector $e$ has at most $.914 \sqrt{m/2}$ nonzero
entries. This is of little practical use here since we are interested
in procedures that might recover a signal when a constant fraction of
the output is unreliable.

Using very different ideas and together with Romberg \cite{CRT}, the
authors proved that the equivalence holds with overwhelming
probability for various types of random matrices provided that
provided that the number of nonzero entries in the vector $e$ be of
the order of $m/\log m$ \cite{RandomBP,OptimalRecovery}. In the
special case where $F$ is an $m/2$ by $m$ random matrix with
independent standard normal entries, \cite{Equivl0l1} proved that the
number of nonzero entries may be as large as $\rho \cdot m$, where
$\rho > 0$ is some very small and unspecified positive constant
independent of $m$.

\subsection{Innovations}

This paper introduces the concept of a \emph{restrictedly almost
  orthonormal system}---a collection of vectors which behaves like an
almost orthonormal system but only for sparse linear combinations.
Thinking about these vectors as the columns of the matrix $F$, we show
that this condition allows for the exact reconstruction of sparse
linear combination of these vectors, i.e. $e$. Our results are
significantly different than those mentioned above as they are
deterministic and do not involve any kind of randomization, although
they can of course be specialized to random matrices.  For instance,
we shall see that a Gaussian matrix with independent entries sampled
from the standard normal distribution is restrictedly almost
orthonormal with overwhelming probability, and that minimizing the
$\ell_1$-norm recovers sparse decompositions with a number of nonzero
entries of size $\rho_0 \cdot m$; we shall actually give numerical
values for $\rho_0$.

We presented the connection with sparse solutions to underdetermined
systems of linear equations merely for pedagogical reasons. There is a
more direct approach. To recover $f$ from corrupted data $y = Af + e$,
we consider solving the following $\ell_1$-minimization problem
\begin{equation}
  \label{eq:P1a}
  (P'_1) \quad \quad  \min_{g \in \R^n} \|y - A g\|_{\ell_1}.
\end{equation}
Now $f$ is the unique solution of $(P'_1)$ if and only if $e$ is the
unique solution of $(P_1)$. In other words, $(P_1)$ and $(P'_1)$ are
equivalent programs. To see why these is true, observe on the one hand
that since $y = Af + e$, we may decompose $g$ as $g = f + h$ so that
\[
(P'_1) \quad \Leftrightarrow  \quad \min_{h \in \R^n} 
\|e - A h\|_{\ell_1}.
\] 
On the other hand, the constraint $F x = F e$ means that $x = e - A h$
for some $h \in \R^n$ and, therefore,
\begin{align*}
  (P_1) \quad & \Leftrightarrow  \quad \min_{h \in \R^n} 
\|x\|_{\ell_1}, \qquad x = e - A h\\
& \Leftrightarrow  \quad \min_{h \in \R^n} 
\|e - Ah\|_{\ell_1},
\end{align*}
which proves the claim. 

The program $(P'_1)$ may also be re-expressed as an LP---hence the
title of this paper. Indeed, the $\ell_1$-minimization problem is
equivalent to
\begin{equation}
   \label{eq:LP}
   \min  1^T t , \qquad -t \le y - A g \le t, 
\end{equation}
where the optimization variables are $t \in R^m$ and $g \in \R^n$ (as
is standard, the generalized vector inequality $x \le y$ means that
$x_i \le y_i$ for all $i$).  As a result, $(P'_1)$ is an LP with
inequality constraints and can be solved efficiently using standard
optimization algorithms, see \cite{BoydBook}.


\subsection{Restricted isometries}

In the remainder of this paper, it will be convenient to use some
linear algebra notations. We denote by $(v_j)_{j \in J} \in R^p$ the
columns of the matrix $F$ and by $H$ the Hilbert space spanned by
these vectors. Further, for any $T \subseteq J$, we let $F_T$ be the
submatrix with column indices $j \in T$ so that
$$
F_T \, c = \sum_{j \in T} c_j v_j \in H. 
$$
To introduce the notion of almost orthonormal system, we first
observe that if the columns of $F$ are sufficiently ``degenerate,''
the recovery problem cannot be solved. In particular, if there exists
a non-trivial sparse linear combination $\sum_{j \in T} c_j v_j = 0$
of the $v_j$ which sums to zero, and $T = T_1 \cup T_2$ is any
partition of $T$ into two disjoint sets, then the vector $y$
$$
y := \sum_{j \in T_1} c_j v_j = \sum_{j \in T_2} (-c_j) v_j
$$
has two distinct sparse representations.  
On the other hand, linear dependencies $\sum_{j \in J} c_j v_j = 0$
which involve a large number of nonzero coefficients $c_j$, as
opposed to a sparse set of coefficients, do not present an obvious
obstruction to sparse recovery. 
At the other extreme, if the $(v_j)_{j \in J}$ are an orthonormal
system, then the recovery problem is easily solved by setting $c_j =
\langle f, v_j \rangle_H$.


The main result of this paper is that if we impose a ``restricted
orthonormality hypothesis,'' which is far weaker than assuming
orthonormality, then $(P_1)$ solves the recovery problem, even if the
$(v_j)_{j \in J}$ are highly linearly dependent (for instance, it is
possible for $m := |J|$ to be much larger than the dimension of
the span of the $v_j$'s).  To make this quantitative we introduce the
following definition.
\begin{definition}[Restricted isometry constants]  
  Let $F$ be the matrix with the finite collection of vectors
  $(v_j)_{j \in J} \in \R^p$ as columns. For every integer $1 \leq S
  \leq |J|$, we define the \emph{$S$-restricted isometry constants}
  $\delta_S$ to be the smallest quantity such that $F_T$ obeys
\begin{equation}\label{frame}
(1-\delta_S) \|c\|^2 \leq \|F_T \, c\|^2 \leq (1 + \delta_S) \|c\|^2
 \end{equation}
 for all subsets $T \subset J$ of cardinality at most $S$, and all
 real coefficients $(c_j)_{j \in T}$.  Similarly, we define the
 \emph{$S,S'$-restricted orthogonality constants} $\theta_{S,S'}$ for
 $S+S' \leq |J|$ to be the smallest quantity such that 
\begin{equation}\label{ortho}
 |\langle F_T c, F_{T'} c' \rangle|  \leq \theta_{S,S'} \cdot  
\|c\| \, \|c'\| 
 \end{equation}
 holds for all \emph{disjoint} sets $T, T' \subseteq J$ of cardinality
 $|T| \leq S$ and $|T'| \leq S'$.
\end{definition}

The numbers $\delta_S$ and $\theta_S$ measure how close the vectors
$v_j$ are to behaving like an orthonormal system, but only when
restricting attention to sparse linear combinations involving no more
than $S$ vectors.  These numbers are clearly non-decreasing in $S$,
$S'$.  For $S = 1$, the value $\delta_1$ only conveys magnitude
information about the vectors $v_j$; indeed $\delta_1$ is the best
constant such that
\begin{equation}\label{length}
 1 - \delta_1 \leq \|v_j \|_H^2 \leq 1 + \delta_1 \hbox{ for all } j \in J.
\end{equation}
In particular, $\delta_1 = 0$ if and only if all of the $v_j$ have
unit length.  Section \ref{sec:uop} establishes that the higher
$\delta_S$ control the orthogonality numbers $\theta_{S,S'}$:
\begin{lemma}\label{theta-delta}  
  We have $\theta_{S,S'} \leq \delta_{S+S'} \leq \theta_{S,S'} +
  \max(\delta_S,\delta_{S'})$ for all $S$, $S'$.
\end{lemma}

To see the relevance of the restricted isometry numbers $\delta_S$ to
the sparse recovery problem, consider the following simple
observation:
\begin{lemma}\label{recover}  
  Suppose that $S \geq 1$ is such that $\delta_{2S} < 1$, and let $T
  \subset J$ be such that $|T| \leq S$.  Let $f := F_T \, c$ for some
  arbitrary $|T|$-dimensional vector $c$.  Then the set $T$ and the
  coefficients $(c_j)_{j \in T}$ can be reconstructed uniquely from
  knowledge of the vector $f$ and the $v_j$'s.
\end{lemma}
\begin{proof}  
  We prove that there is a unique $c$ with $\|c\|_{\ell_0} \le S$ and
  obeying $f = \sum_j c_j v_j$.  Suppose for contradiction that $f$
  had two distinct sparse representations $f = F_T \, c = F_{T'} \,
  c'$ where $|T|, |T'| \leq S$. Then
$$
F_{T \cup T'}\, d = 0, \qquad d_j := c_j 1_{j \in T} - c_{j'} 1_{j \in T'}.
$$
Taking norms of both sides and applying \eqref{frame} and the
hypothesis $\delta_{2S} < 1$ we conclude that $\|d\|^2 = 0$,
contradicting the hypothesis that the two representations were
distinct.
\end{proof}

\subsection{Main results}

Note that the previous lemma is an abstract existence argument which
shows what might theoretically be possible, but does not supply any
efficient algorithm to recover $T$ and $c_j$ from $f$ and $(v_j)_{j
  \in J}$ other than by brute force search---as discussed earlier.  In
contrast, our main theorem result that, by imposing slightly stronger
conditions on $\delta_{2S}$, the $\ell_1$-minimization program $(P_1)$
recovers $f$ exactly.
\begin{theorem}\label{recover-better} 
Suppose that $S \geq 1$ is such that 
\begin{equation}\label{delta-cond}
\delta_S + \theta_S + \theta_{S,2S} < 1,
\end{equation} 
and let $c$ be a real vector supported on a set $T \subset J$ obeying
$|T| \leq S$.  Put $f := F c$.  Then $c$ is the unique minimizer to
\[
(P_1) \quad \quad \min \|d\|_{\ell_1} \qquad F d = f. 
\]
\end{theorem}
Note from Lemma \ref{theta-delta} that \eqref{delta-cond} implies
$\delta_{2S} < 1$, and is in turn implied by $\delta_S +
\delta_{2S}+\delta_{3S} < 1/4$.  Thus the condition \eqref{delta-cond}
is roughly ``three times as strict'' as the condition required for
Lemma \ref{recover}.

Theorem \ref{recover-better} is inspired by our previous work
\cite{OptimalRecovery}, see also \cite{CRT,RandomBP}, but unlike those
earlier results, our results here are \emph{deterministic}, and thus
do not have a non-zero probability of failure, provided of course one
can ensure that the system $(v_j)_{j \in J}$ verifies the condition
\eqref{delta-cond}. By virtue of the previous discussion, we have the
companion result:
\begin{theorem}
\label{teo:decode} 
Suppose $F$ is such that $F A = 0$ and let $S \geq 1$ be a number
obeying the hypothesis of Theorem \ref{recover-better}. Set $y = A f +
e$, where $e$ is a real vector supported on a set of size at most $S$.
Then $f$ is the unique minimizer to
\[
(P'_1) \quad \quad \min_{g \in \R^n} \|y - A g\|_{\ell_1}.  
\]
\end{theorem}

\subsection{Gaussian random matrices}

An important question is then to find matrices with good restricted
isometry constants, i.e. such that \eqref{delta-cond} holds for large
values of $S$. Indeed, such matrices will tolerate a larger fraction
of output in error while still allowing exact recovery of the original
input by linear programming. How to construct such matrices might be
delicate. In section \ref{sec:random}, however, we will argue that
generic matrices, namely samples from the Gaussian unitary ensemble
obey $\eqref{delta-cond}$ for relatively large values of $S$.
\begin{theorem} 
\label{teo:gauss}
Assume $p \le m$ and let $F$ be a $p$ by $m$ matrix whose entries are
i.i.d. Gaussian with mean zero and variance $1/p$. Then the condition
of Theorem \ref{recover-better} holds with overwhelming probability
provided that $r = S/m$ is small enough so that
\[
r < r^*(p,m)
\]
where $r^*(p,m)$ is given in Section \ref{eq:rho}. (By ``with
overwhelming probability,'' we mean with probability decaying
exponentially in $m$.)  In the limit of large samples, $r^*$ only
depends upon the ratio, and numerical evaluations show that the
condition holds for $r \le 3.6 \cdot 10^{-4}$ in the case where $p/m =
3/4$, $r \le 3.2 \cdot 10^{-4}$ when $p/m = 2/3$, and $r \le 2.3 \cdot
10^{-4}$ when $p/m = 1/2$.
\end{theorem}
In other words, Gaussian matrices are a class of matrices for which
one can solve an underdetermined systems of linear equations by
minimizing $\ell_1$ provided, of course, the input vector has fewer
than $\rho \cdot m$ nonzero entries with $\rho > 0$. We mentioned
earlier that this result is similar to \cite{Equivl0l1} . What is new
here is that by using a very different machinery, one can obtain
explicit numerical values which were not available before.

In the context of error correcting, the consequence is that a fraction
of the output may be corrupted by {\em arbitrary} errors and yet,
solving a convex problem would still recover $f$ exactly---a rather
unexpected feat.
\begin{corollary}
  Suppose $A$ is an $n$ by $m$ Gaussian matrix and set $p = m - n$.
  Under the hypotheses of Theorem \ref{teo:gauss}, the solution to
  $(P_1')$ is unique and equal to $f$. 
\end{corollary}
This is an immediate consequence of Theorem \ref{teo:gauss}. The only
thing we need to argue is why we may think of the annihilator $F$
(such that $FA = 0$) as a matrix with independent Gaussian entries.
Observe that the range of $A$ is a random space of dimension $n$
embedded in $\R^m$ so that the data $\tilde y = F e$ is the projection
of $e$ on a random space of dimension $p$.  The range of a $p$ by $m$
matrix with independent Gaussian entries precisely is a random
subspace of dimension $p$, which justifies the claim.

We would like to point out that the numerical bounds we derived in
this paper are overly pessimistic. We are confident that finer
arguments and perhaps new ideas will allow to derive versions of
Theorem \ref{teo:gauss} with better bounds. The discussion section
will enumerate several possibilities for improvement.

\subsection{Organization of the paper}

The paper is organized as follows. Section \ref{sec:main} proves our
main claim, namely, Theorem \ref{recover-better} (and hence Theorem
\ref{teo:decode}) while Section \ref{sec:random} introduces elements
from random matrix theory to establish Theorem \ref{teo:gauss}. In
Section \ref{sec:numerical}, we present numerical experiments which
suggest that in practice, $(P'_1)$ works unreasonably well and
recovers the $f$ exactly from $y = Af + e$ provided that the fraction
of the corrupted entries be less than about 17\% in the case where $m
= 2n$ and less than about 34\% in the case where $m = 4n$. Section
\ref{sec:optimalrecovery} explores the consequences of our results for
the recovery of signals from highly incomplete data and ties our
findings with some of our earlier work. Finally, we conclude with a
short discussion section whose main purpose is to outline areas for
improvement.

\section{Proof of Main Results}
\label{sec:main}

Our main result, namely, Theorem \ref{recover-better} is proved by
duality. As we will see in section \ref{sec:main-proof}, $c$ is the
unique minimizer if the matrix $F_T$ has full rank and if one can find
a vector $w$ with the two properties
\begin{itemize}
\item[(i)] $\langle w, v_j \rangle_H = \sgn(c_j)$ for all $j \in T$, 
\item[(ii)] and $|\langle w, v_j \rangle_H| < 1$ for all $j
  \notin T$,
\end{itemize}
where $\sgn(c_j)$ is the sign of $c_j$ ($\sgn(c_j) = 0$ for $c_j =
0$). The two conditions above say that a specific dual program is
feasible and is called the {\em exact reconstruction property} in
\cite{OptimalRecovery}, see also \cite{CRT}.  For $|T| \le S$ with $S$
obeying the hypothesis of Theorem \ref{recover-better}, $F_T$ has full
rank since $\delta_S < 1$ and thus, the proof simply consists in
constructing a dual vector $w$; this is the object of the next
section.

\subsection{Exact reconstruction property}

We now examine the sparse reconstruction property and begin with
coefficients $\langle w, v_j \rangle_H$ for $j \not \in T$ being only
small in an $\ell_2$ sense.

\begin{lemma}[Dual sparse reconstruction property, $\ell_2$ version]
\label{erp-2}  
Let $S,S' \geq 1$ be such that $\delta_S < 1$, and $c$ be a real
vector supported on $T \subset J$ such that $|T| \leq S$. Then there
exists a vector $w \in H$ such that $\langle w, v_j \rangle_H = c_j$
for all $j \in T$.  Furthermore, there is an ``exceptional set''
$E \subset J$ which is disjoint from $T$, of size at most 
\begin{equation}\label{et}
|E| \leq S', 
\end{equation}
and with the properties 
$$
|\langle w, v_j \rangle| \leq
\frac{\theta_{S,S'}}{(1-\delta_S)\sqrt{S'}} \cdot \|c\| \text{ for
  all } j \not \in T \cup E$$
and
$$
(\sum_{j \in E} |\langle w, v_j \rangle|^2)^{1/2} \leq
\frac{\theta_S}{1 - \delta_{S}} \cdot \|c\|.$$
In addition, $\|w\|_H
\leq K \cdot \|c\|$ for some constant $K > 0$ only depending upon
$\delta_S$. 
\end{lemma}
\begin{proof}  Recall that $F_T: \ell_2(T) \to H$ is the linear 
  transformation $F_T \, c_T := \sum_{j \in T} c_j v_j$ where $c_T :=
  (c_j)_{j \in T}$ (we use the subscript $T$ in $c_T$ to emphasize
  that the input is a $|T|$-dimensional vector), and let $F_T^*$ be the
  adjoint transformation
$$
  F_T^* w := (\langle w, v_j \rangle_H)_{j \in T}.
$$
  Property 
  \eqref{frame} gives   
  $$
  1 - \delta_{S} \leq \lambda_{\min}(F_T^* F_T) \leq
  \lambda_{\max}(F_T^* F_T) \leq 1 + \delta_{S},$$
  where
  $\lambda_{\min}$ and $\lambda_{\max}$ are the minimum and maximum
  eigenvalues of the positive-definite operator $F_T^* F_T$.  In
  particular, since $\delta_{|T|} < 1$, we see that $F_T^* F_T$ is
  invertible with
\begin{equation}\label{fft-inv}
 \| (F^*_T F_T)^{-1} \| \leq \frac{1}{1 - \delta_{S}}.
 \end{equation}
 Also note that $\| F_T (F^*_T F_T)^{-1}\| \le \sqrt{1 +
   \delta_S}/(1-\delta_S)$ and set $w \in H$ to be the vector
 $$
 w := F_T (F^*_T F_T)^{-1} c_T;
 $$
 it is then clear that $F^*_T w = c_T$, i.e.  $\langle w, v_j
 \rangle_H = c_j$ for all $j \in T$.  In addition, $\|w\| \le K \cdot
 \|c_T\|$ with $K = \sqrt{1 + \delta_S}/(1-\delta_S)$. 
 Finally, if $T'$ is any set in $J$ disjoint from $T$ with $|T'| \leq
 S'$ and $d_{T'} = (d_j)_{j \in T'}$ is any sequence of real numbers,
 then \eqref{ortho} and \eqref{fft-inv} give
\begin{align*}
  |\langle F^*_{T'} w, d_{T'} \rangle_{\ell_2(T')}| = |\langle w,
  F_{T'} d_{T'} \rangle_{\ell_2(T')}| & = \left|\langle \sum_{j \in T}
    ((F^*_T F_T)^{-1} c_T)_j v_j,
    \sum_{j \in T'} d_j v_j \rangle_H\right| \\
  &\leq \theta_{S,S'} \cdot \| (F^*_T F_T)^{-1} c_T \| \cdot
  \| d_{T'} \| \\
  &\leq \frac{\theta_{S,S'}}{1 - \delta_S} \| c_T \| \cdot \|d_{T'}\|;
\end{align*}
since $d_{T'}$ was arbitrary, we thus see from duality that
$$
\| F^*_{T'} w \|_{\ell_2(T')} \leq \frac{\theta_{S,S'}}{1 -
  \delta_{S}} \| c_T\|.$$
In other words,
\begin{equation}\label{tbone}
 (\sum_{j \in T'} |\langle w, v_j \rangle|^2)^{1/2} \leq
\frac{\theta_{S,S'}}{1 - \delta_{S}} \| c_T \| 
\hbox{ whenever } T' \subset J \backslash T \hbox{ and } |T'| \leq S'.
\end{equation}
If in particular if we set
$$
E := \{ j \in J \backslash T: |\langle w, v_j \rangle| >
\frac{\theta_{S,S'}}{(1 - \delta_{S})\sqrt{S'}} \cdot \|c_T\|\},
$$
then $|E|$ must obey $|E| \leq S'$, since otherwise we could
contradict \eqref{tbone} by taking a subset $T'$ of $E$ of cardinality
$S'$.  The claims now follow.
\end{proof}

We now derive a solution with better control on the sup norm of $|\<
w, v_j\>|$ outside of $T$, by iterating away the exceptional set $E$
(while keeping the values on $T$ fixed).

\begin{lemma}[Dual sparse reconstruction property, $\ell_\infty$ version]
\label{erp-3} 
Let $S \geq 1$ be such that $\delta_S + \theta_{S,2S} < 1$, and $c$ be
a real vector supported on $T \subset J$ obeying $|T| \le S$. Then
there exists a vector $w \in H$ such that $\langle w, v_j \rangle_H =
c_j$ for all $j \in T$.  Furthermore, $w$ obeys
\begin{equation}
  \label{eq:erp-3}
|\langle w, v_j \rangle| \leq \frac{\theta_S}{(1 - \delta_S -
  \theta_{S,2S}) \sqrt{S}} \cdot \|c\| \text{ for all } j \not \in
T.
\end{equation}
\end{lemma}
\begin{proof}  
  We may normalize $\sum_{j \in T} |c_j|^2 = \sqrt{S}$.  Write $T_0 :=
  T$.  Using Lemma \ref{erp-2}, we can find a vector $w_1 \in H$ and a
  set $T_1 \subseteq J$ such that
\begin{align*}
  T_0 \cap T_1 &= \emptyset\\
  |T_1| &\leq S \\
  \langle w_1, v_j \rangle_H &= c_j \hbox{ for all } j \in T_0\\
  |\langle w_1, v_j \rangle_H| &\leq \frac{\theta_{S,S'}}{(1-\delta_S)} 
\hbox{ for all } j \not \in T_0 \cup T_1 \\
  (\sum_{j \in T_1} |\langle w_1, v_j \rangle_H|^2)^{1/2} 
&\leq \frac{\theta_S}{1 - \delta_{S}} \sqrt{S}\\
  \|w_1\|_H &\leq K.
\end{align*}
Applying Lemma \ref{erp-2} iteratively gives a sequence of vectors
$w_{n+1} \in H$ and sets $T_{n+1} \subseteq J$ for all $n \geq 1$ with
the properties 
\begin{align*}
  T_n \cap (T_0 \cup T_{n+1}) &= \emptyset\\
  |T_{n+1}| &\leq S \\
  \langle w_{n+1}, v_j \rangle_H &= \langle w_n, v_j \rangle_H
  \hbox{ for all } j \in T_n\\
  \langle w_{n+1}, v_j \rangle_H &= 0 \hbox{ for all } j \in T_0\\
  |\langle w_{n+1}, v_j \rangle_H| &\leq \frac{\theta_S}{1 -
    \delta_{S}} \left(\frac{\theta_{S,2S}}{1-\delta_S}\right)^{n}
  \hbox{ for all } j \not \in T_0 \cup T_n \cup T_{n+1} \\
  (\sum_{j \in T_{n+1}} |\langle w_{n+1}, v_j \rangle|^2)^{1/2} &\leq
  \frac{\theta_S}{1 - \delta_{S}}
  \left(\frac{\theta_{S,2S}}{1 - \delta_{S}}\right)^{n} \sqrt{S}\\
  \| w_{n+1} \|_H &\leq \left(\frac{\theta_S}{1 -
      \delta_{S}}\right)^{n-1} K.
\end{align*}
By hypothesis, we have $\frac{\theta_{S,2S}}{1-\delta_S} \leq 1$.
Thus if we set
\[
w := \sum_{n=1}^\infty (-1)^{n-1} w_n
\]
then the series is absolutely convergent and, therefore, $w$ is a
well-defined vector in $H$.  We now study the coefficients
\begin{equation}\label{w-converge}
\langle w, v_j \rangle_H = 
\sum_{n=1}^\infty (-1)^{n-1} \langle w_n, v_j\rangle_H
\end{equation}
for $j \in J$.

Consider first $j \in T_0$, it follows from the construction that
$\langle w_1, v_j \rangle_H = c_j$ and $\langle w_n, v_j \rangle_H =
0$ for all $n \geq 2$, and hence
\[
\langle w, v_j \rangle_H = c_j \hbox{ for all } j \in T_0.
\]
Second, fix $j$ with $j \not \in T_0$ and let $I_j := \{ n \geq 1: j
\in T_n \}$.  Since $T_n$ and $T_{n+1}$ are disjoint, we see that the
integers in the set $I_j$ are spaced at least two apart.  Now if $n
\in I_j$, then by definition $j \in T_n$ and, therefore,
\[
\langle w_{n+1}, v_j \rangle_H = \langle w_n, v_j \rangle_H.
\]
In other words, the $n$ and $n+1$ terms in \eqref{w-converge} cancel
each other out.  Thus we have
\[
\langle w, v_j \rangle_H = \sum_{n \geq 1; n, n-1 \not \in I_j}
(-1)^{n-1} \langle w_n, v_j\rangle_H.
\]
On the other hand, if $n, n-1
\not \in I_j$ and $n \neq 0$, then $j \not \in T_n \cap T_{n-1}$ and 
$$
|\langle w_{n}, v_j \rangle| \leq \frac{\theta_{S,S}}{1 -
  \delta_{S}} \left(\frac{\theta_{S,2S}}{1-\delta_S}\right)^{n-1}$$
which by the triangle inequality and the geometric series formula
gives
$$
|\sum_{n \geq 1; n, n-1 \not \in I_j} (-1)^{n-1} \langle w_n,
v_j\rangle_H| \leq \frac{\theta_{S,S}}{1 - \delta_S -
  \theta_{S,2S}}.$$
In conclusion, 
$$
|\langle w, v_j \rangle_H - 1_{0 \in I_j} \langle w_0, v_j
\rangle_H| \leq \frac{\theta_{S,S}}{1 - \delta_S - \theta_{S,2S}},$$
and since $|T| \leq S$, the claim follows.
\end{proof}

Lemma \ref{erp-3} actually solves the dual recovery problem. Indeed,
our result states that one can find a vector $w \in H$ obeying both
properties (i) and (ii) stated at the beginning of the section. To see
why (ii) holds, observe that $\|\sgn(c)\| = \sqrt{|T|} \le \sqrt{S}$
and, therefore, \eqref{eq:erp-3} gives for all $j \not \in T$
\begin{equation*}
  |\langle w, v_j \rangle_H| \leq 
\frac{\theta_{S,S}}{(1 - \delta_S - \theta_{S,2S})} \cdot 
\sqrt{\frac{|T|}{S}} 
  \leq \frac{\theta_{S,S}}{(1 - \delta_S - \theta_{S,2S})} < 1,
\end{equation*}
provided that $\delta_S + \theta_{S,S} + \theta_{S,2S} < 1$.

\subsection{Proof of Theorem \ref{recover-better}}  
\label{sec:main-proof}

Observe first that standard convex arguments give that there exists at
least one minimizer $d = (d_j)_{j \in J}$ to the problem $(P_1)$.  We
need to prove that $d = c$.  Since $c$ obeys the constraints of this
problem, $d$ obeys
\begin{equation}\label{d-meander}
 \| d \|_{\ell_1} \leq 
\|c\|_{\ell_1} = \sum_{j \in T} |c_j|.
 \end{equation}
 Now take a $w$ obeying properties (i) and (ii) (see the remark
 following Lemma \ref{erp-3}).  Using the fact that the inner product
 $\langle w, v_j\rangle$ is equal to the sign of $c$ on $T$ and has
 absolute value strictly less than one on the complement, we then
 compute
\begin{align*}
  \| d \|_{\ell_1} &= 
\sum_{j \in T} |c_j + (d_j-c_j)| + \sum_{j \not \in T} |d_j| \\
  &\geq \sum_{j \in T} \sgn(c_j) (c_j + (d_j-c_j)) + 
\sum_{j \not \in T} d_j \langle w, v_j \rangle_H \\
  &= \sum_{j \in T} |c_j| + \sum_{j \in T} (d_j - c_j) 
\langle w, v_j \rangle_H + \sum_{j \not \in T} d_j \langle w, v_j \rangle_H  \\
  &= \sum_{j \in T} |c_j| + 
\langle w, \sum_{j \in J} d_j v_j - \sum_{j \in T} c_j \rangle \\
  &= \sum_{j \in T} |c_j| + \langle w, f - f \rangle \\
  &= \sum_{j \in T} |c_j|.
\end{align*}
Comparing this with \eqref{d-meander} we see that all the inequalities
in the above computation must in fact be equality.  Since $|\langle w,
v_j \rangle_H|$ was strictly less than 1 for all $j \not \in T$, this
in particular forces $d_j = 0$ for all $j \notin T$.  Thus
$$
\sum_{j \in T} (d_j - c_j) v_j = f - f = 0.$$
Applying
\eqref{frame} (and noting from hypothesis that $\delta_S < 1$) we
conclude that $d_j = c_j$ for all $j \in T$.  Thus $d = c$ as claimed.
This concludes the proof of our theorem.

{\bf Remark.}  It is likely that one may push the condition $\delta_S
+ \theta_{S,S} + \theta_{S,2S} < 1$ a little further.  The key idea is
as follows. Each vector $w_n$ in the iteration scheme used to prove
Lemma \ref{erp-3} was designed to annihilate the influence of
$w_{n-1}$ on the exceptional set $T_{n-1}$.  But annihilation is too
strong of a goal.  It would be just as suitable to design $w_n$ to
moderate the influence of $w_{n-1}$ enough so that the inner product
with elements in $T_{n-1}$ is small rather than zero. However, we have 
not pursued such refinements as the arguments would become considerably 
more complicated than the calculations presented here.

\subsection{Approximate orthogonality}
\label{sec:uop}

Lemma \ref{theta-delta} gives control of the size of the principal
angle between subspaces of dimension $S$ and $S'$ respectively. This
is useful because it allows to guarantee exact reconstruction from the
knowledge of the $\delta$ numbers only.

\begin{proof}[Proof of Lemma \ref{theta-delta}]
  We first show that $\theta_{S,S'} \leq \delta_{S+S'}$.  By
  homogeneity it will suffice to show that
  $$
  |\langle \sum_{j \in T} c_j v_j, \sum_{j' \in T'} c'_{j'} v_{j'}
  \rangle_{H}| \leq \delta_{S+S'}$$
  whenever $|T| \leq S$, $|T'| \leq
  S'$, $T, T'$ are disjoint, and $\sum_{j \in T} |c_j|^2 = \sum_{j'
    \in T'} |c'_{j'}|^2 = 1$. Now \eqref{frame} gives 
  $$
  2(1-\delta_{S+S'}) \leq \| \sum_{j \in T} c_j v_j + \sum_{j' \in
    T'} c'_{j'} v_{j'} \|_H^2 \leq 2(1+\delta_{S+S'})$$
  together with 
  $$
  2(1-\delta_{S+S'}) \leq \| \sum_{j \in T} c_j v_j - \sum_{j' \in
    T'} c'_{j'} v_{j'} \|_H^2 \leq 2(1+\delta_{S+S'}),$$
  and the claim now
  follows from the parallelogram identity
  $$
  \langle f, g \rangle = \frac{\|f+g\|_H^2 - \|f-g\|_H^2}{4}.
  $$
  It remains to show that $\delta_{S+S'} \leq \theta_S + \delta_S$.
  Again by homogeneity, it suffices to establish that
  $$
  |\langle \sum_{j \in \tilde T} c_j v_j, \sum_{j' \in \tilde T}
  c_{j'} v_{j'} \rangle_{H} - 1| \leq (\delta_{S} + \theta_S)$$
  whenever $|\tilde T| \leq S+S'$ and $\sum_{j \in \tilde T} |c_j|^2 =
  1$.  To prove this property, we partition $\tilde T$ as $\tilde T =
  T \cup T'$ where $|T| \leq S$ and $|T'| \leq S'$ and write $\sum_{j
    \in T} |c_j|^2 = \alpha$ and $\sum_{j \in T'} |c_j|^2 = 1-\alpha$.
  \eqref{frame} together with \eqref{ortho} give
\begin{align*}
  (1-\delta_S) \alpha \leq \langle \sum_{j \in T} c_j v_j, 
\sum_{j' \in T} c_{j'} v_{j'} \rangle_{H} &\leq (1+\delta_S) \alpha,\\
  (1-\delta_{S'}) (1-\alpha) \leq \langle \sum_{j \in T'} c_j v_j, 
\sum_{j' \in T'} c_{j'} v_{j'} \rangle_{H} &\leq (1+\delta_{S'}) (1-\alpha),\\
  |\langle \sum_{j \in T} c_j v_j, \sum_{j' \in T} c_{j'} v_{j'}
  \rangle_{H}| &\leq \theta_{S,S'} \alpha^{1/2} (1-\alpha)^{1/2}. 
\end{align*}
Hence
\begin{align*}
  |\langle \sum_{j \in \tilde T} c_j v_j, \sum_{j' \in \tilde T}
  c_{j'} v_{j'} \rangle_{H} - 1| &\leq \delta_S \alpha + \delta_{S'}
  (1-\alpha)
  + 2 \theta_S \alpha^{1/2} (1-\alpha)^{1/2} \\
  &\leq \max(\delta_S,\delta_{S'}) + \theta_S
\end{align*}
as claimed.  (We note that it is possible to optimize this bound a
little further but will not do so here.)
\end{proof}

\section{Gaussian Random Matrices}
\label{sec:random}

In this section, we argue that with overwhelming probability, Gaussian
random matrices have ``good'' isometry constants. Consider a
$p$ by $m$ matrix $F$ whose entries are i.i.d. Gaussian with mean zero
and variance $1/p$ and let $T$ be a subset of the columns. We wish to
study the extremal eigenvalues of $F_T^* F_T$.
Following upon the work of Marchenko and Pastur
\cite{MarchenkoPastur}, Geman \cite{Geman} and Silverstein
\cite{Silverstein} (see also \cite{BaiYin}) proved that
\begin{align*}
  \label{eq:geman}
 & \lambda_{\min}(F_T^* F_T) \goto (1 - \sqrt{\gamma})^2 
\text{ a.s.}\\
 & \lambda_{\max}(F_T^* F_T) \goto (1 + \sqrt{\gamma})^2  \text{ a.s.},
\end{align*}
in the limit where $p$ and $|T| \goto \infty$ with
\[
|T|/p \goto \gamma \le 1.
\]
In other words, this says that loosely speaking and in the limit of
large $p$, the restricted isometry constant $\delta(F_T)$ for a fixed
$T$ behaves like
\[
1 - \delta(F_T) \le \lambda_{\min}(F_T^* F_T) \le \lambda_{\max}(F_T)
\le 1 + \delta(F_T), \qquad \delta(F_T) \approx 2 \sqrt{|T|/p} + |T|/p.
\]

Restricted isometry constants must hold for all sets $T$ of
cardinality less or equal to $S$, and we shall make use of
concentration inequalities to develop such a uniform bound. Note that
for $T' \subset T$, we obviously have
\[
\lambda_{\min}(F_T^* F_T) \le \lambda_{\min}(F_{T'}^* F_{T'}) \le
\lambda_{\max}(F_{T'}^* F_{T'}) \le \lambda_{\max}(F_T^* F_T)
\]
and, therefore, attention may be restricted to matrices of size $S$.
Now, there are large deviation results about the singular values of
$F_T$ \cite{Szarek2}. For example, letting $\sigma_{\max}(F_T)$
(resp.~$\sigma_{\min}$) be the largest singular value of $F_T$ so that
$\sigma^2_{\max}(F_T) = \lambda_{\max}(F_T^*F_T)$
(resp.~$\sigma^2_{\min}(F_T) = \lambda_{\min}(F_T^*F_T)$), Ledoux
\cite{LedouxAMS} applies the concentration inequality for Gaussian
measures, and for a each fixed $t > 0$, obtains the deviation bounds
\begin{align}
\label{eq:concentration}
&\P\left(\sigma_{\max}(F_T)   >   
1 + \sqrt{|T|/p} + o(1) + t\right) \le e^{- p t^2/2}\\
\label{eq:concentration2}
&\P\left(\sigma_{\min}(F_T)   < 
1 - \sqrt{|T|/p} + o(1) - t \right) \le e^{- p t^2/2};
\end{align}
here, $o(1)$ is a small term tending to zero as $p \goto \infty$ and
which can be calculated explicitly, see \cite{ElKarouiPhD}. For
example, this last reference shows that one can select $o(1)$ in
\eqref{eq:concentration} as $\frac{1}{2p^{1/3}} \cdot
\gamma^{1/6}(1+\sqrt{\gamma})^{2/3}$.

\begin{lemma}
\label{teo:Gauss_delta_S}
Put $r = S/m$ and set 
\[
f(r) := \sqrt{m/p} \cdot \left(\sqrt{r} + \sqrt{2H(r)}\right), 
\] 
where $H$ is the entropy function $H(q) := - q\log q - (1-q)
\log(1-q)$ defined for $0 < q < 1$.  For each $\epsilon> 0$, the
restricted isometry constant $\delta_S$ of a $p$ by $m$ Gaussian
matrix $F$ obeys
\begin{equation}
\label{large-dev}
  \P\left(1 + \delta_S > [1 + (1 + \epsilon) f(r)]^2\right) \le 2 
\cdot e^{-m H(r) \cdot \epsilon/2}.
\end{equation}
\end{lemma}
\begin{proof}
  As discussed above, we may restrict our attention to sets $|T|$ such
  that $|T| = S$. Denote by $\eta_p$ the $o(1)$-term appearing in
  either \eqref{eq:concentration} or \eqref{eq:concentration2}. Put
  $\lambda_{\max} = \lambda_{\max}(F_T^* F_T)$ for short, and observe
  that
\begin{eqnarray*}
\P\left(\sup_{T: |T| = S} \lambda_{\max} > 
(1 + \sqrt{S/p} + \eta_p + t)^2\right) & \le &  |\{T: |T| = S\}|  
\, \P\left(\lambda_{\max} > (1 + \sqrt{S/p} + \eta_p + t)^2\right) \\
& \le &  \binom{m}{S} \, e^{-pt^2/2}.  
\end{eqnarray*}
>From Stirling's approximation $\log m! = m \log m - m + O(\log m)$ we have
the well-known formula
\[
\log \binom{m}{S} = m H(r) + O(\log m).
\]
which gives
\[
P\left(\sup_{T: |T| = S} \lambda_{\max} > 
(1 + \sqrt{S/p} +  \eta_p + t)^2\right) 
\le e^{m H(r)} \cdot e^{O(\log m)} \cdot e^{-pt^2/2},
\]
The exact same argument
applied to the smallest eigenvalues yields
\[
P\left(\inf_{T: |T| = S} \lambda_{\min} < (1 - \sqrt{S/p} -\eta_p -
  t)^2\right) \le e^{m H(r)} \cdot e^{O(\log(m))} \cdot e^{-pt^2/2}.
\]
Fix $\eta_p + t = (1+\epsilon) \cdot \sqrt{m/p} \cdot \sqrt{2H(r)}$. Assume
now that $m$ and $p$ are large enough so that 
$\eta_p \le \epsilon/2 \cdot \sqrt{m/p} \cdot \sqrt{2H(r)}$. Then
\[
P\left(\sup_{T: |T| = S} \lambda_{\max} > (1 + \sqrt{S/p} + \epsilon
  \cdot \sqrt{m/p} \cdot \sqrt{2H(r)})^2\right) \le e^{-m H(r) \cdot
  \epsilon/2}.
\]
where we used the fact that the term $O(\log m)$ is less than $m
\epsilon H(r)/2$ for sufficiently large $m$.  The same bound holds for
the minimum eigenvalues and the claim follows.
\end{proof}

\begin{figure}
  \centering 
\includegraphics[width=3.5in]{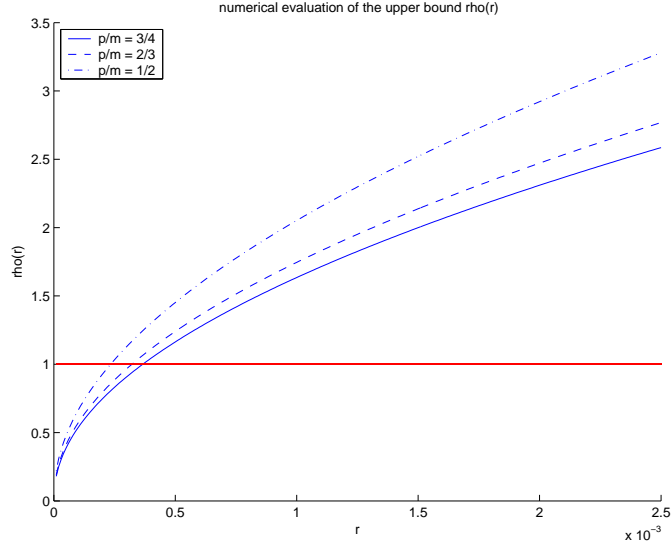}
  \caption{Behavior of the upper bound $\rho_{p/m}(r)$ for three
    values of the ratio $p/m$, namely, $p/m = 3/4, 2/3, 1/2$.}
 \label{fig:bound} 
\end{figure}

Ignoring the $\epsilon$'s, Lemma \ref{teo:Gauss_delta_S} states that
with overwhelming probability
\begin{equation}
  \label{eq:f}
  \delta_S < - 1 + [1 + f(r)]^2.
\end{equation}
A similar conclusion holds for $\delta_{2S}$ and $\delta_{3S}$ and,
therefore, we established that
\begin{equation}
  \label{eq:rho}
\delta_S + \delta_{2S} + \delta_{3S} 
< \rho_{p/m}(r), \qquad \rho_{p/m}(r) =  \sum_{j = 1}^3 -1 + [1 + f(jr)]^2. 
\end{equation}
with very high probability. In conclusion, Lemma \ref{theta-delta}
shows that the hypothesis of our main theorem holds provided that the
ratio $r = S/m$ be small so that $ \rho_{p/m}(r) < 1$. In other words,
in the limit of large samples, $S/m$ maybe taken as any value obeying
$\rho_{p/m}(S/m) < 1$ which we used to give numerical values in Theorem
\ref{teo:gauss}.  Figure~\ref{fig:bound} graphs the function
$\rho_{p/m}(r)$ for several values of the ratio $p/m$.

\section{Numerical Experiments}
\label{sec:numerical}

This section investigates the practical ability of $\ell_1$ to recover
an object $f \in \R^n$ from corrupted data $y = A f + e$, $y \in \R^m$
(or equivalently to recover the sparse vector of errors $e \in \R^m$
from the underdetermined system of equations $F e = z \in \R^{m-n}$).
The goal here is to evaluate empirically the location of the
breakpoint as to get an accurate sense of the performance one might
expect in practice. In order to do this, we performed a series of
experiments designed as follows:
\begin{enumerate}
\item select $n$ (the size of the input signal) and $m$ so that with
  the same notations as before, $A$ is an $n$ by $m$ matrix; sample
  $A$ with independent Gaussian entries; 
    
\item select $S$ as a percentage of $m$; 
  
\item select a support set $T$ of size $|T|= S$ uniformly at random,
  and sample a vector $e$ on $T$ with independent and identically
  distributed Gaussian entries\footnote{Just as in \cite{RandomBP},
    the results presented here do not seem to depend on the actual
    distribution used to sample the errors.};
\item make $\tilde y = A x + e$ (the choice of $x$ does not matter as
  is clear from the discussion and here, $x$ is also selected at
  random), solve $(P'_1)$ and obtain $x^*$;
\item compare $x$ to $x^*$;
\item repeat $100$ times for each $S$ and $A$;
\item repeat for various sizes of $n$ and $m$. 
\end{enumerate}

\begin{figure}
\centerline{
\begin{tabular}{ccc}
\includegraphics[width=3in]{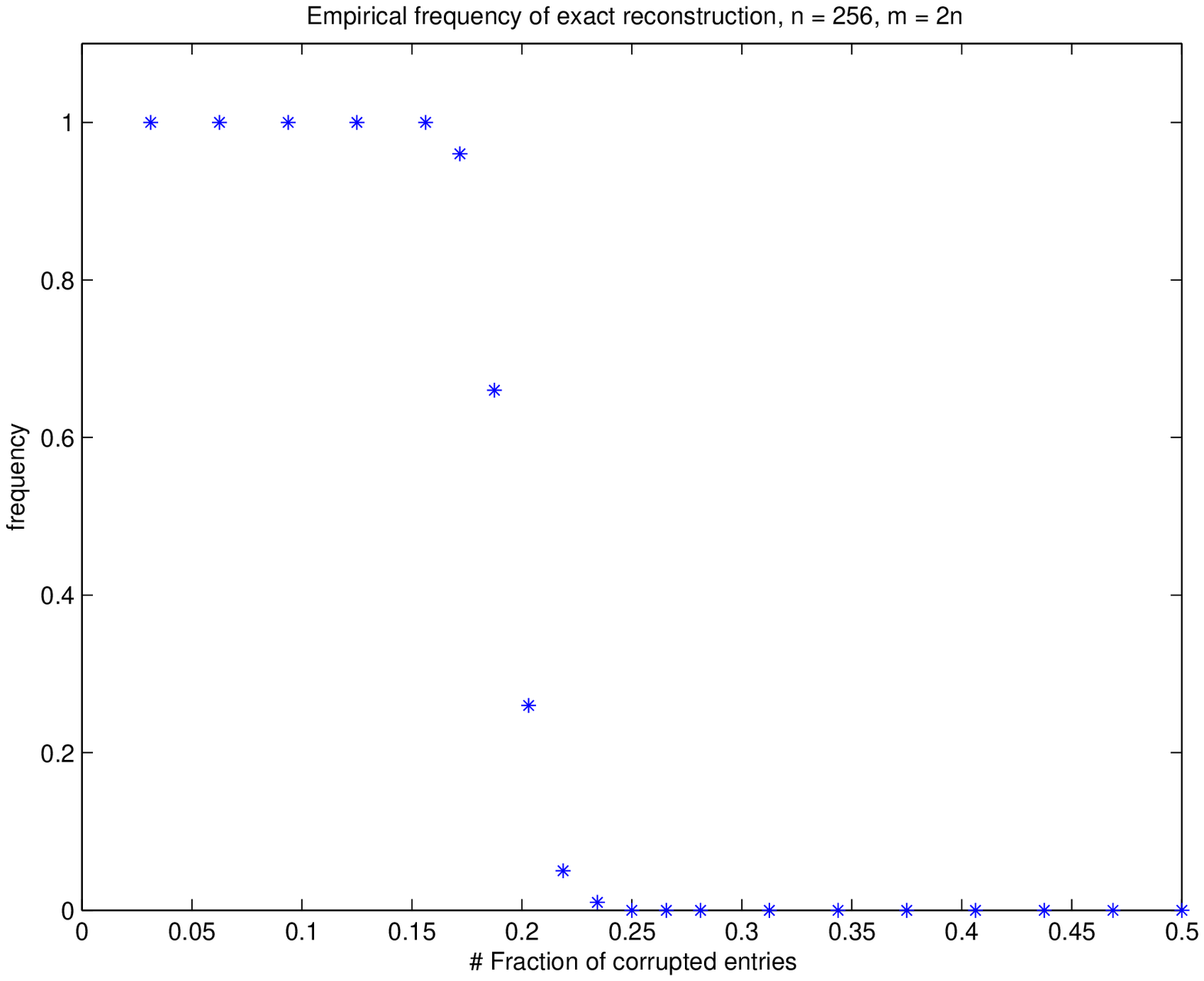} & \hspace{5mm} &
\includegraphics[width=3in]{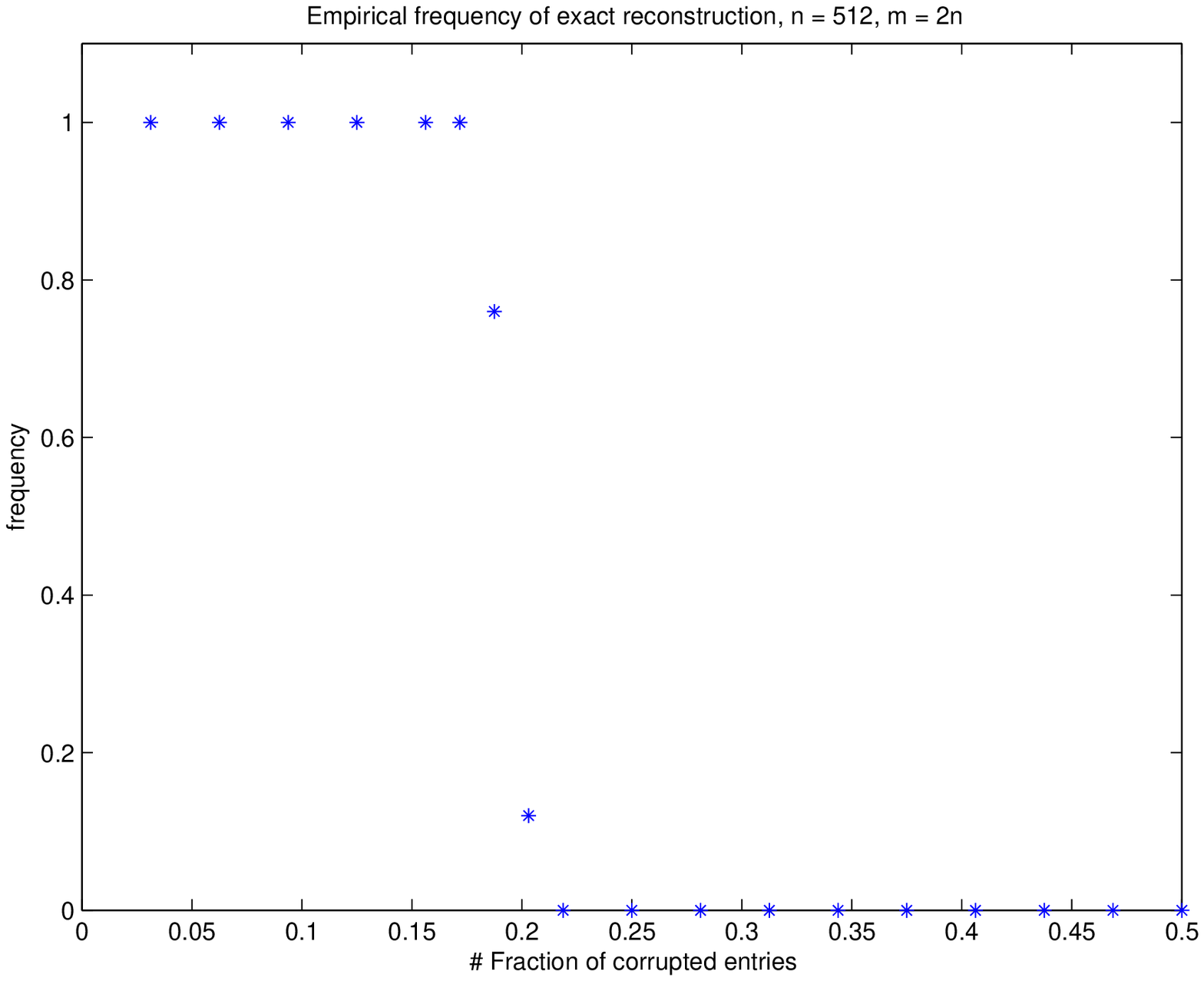} \\
(a) &  & (b)
\end{tabular}
}
\caption{\small\sl $\ell_1$-recovery of an input signal from 
  $y = Af + e$ with $A$ an $m$ by $n$ matrix with independent Gaussian
  entries.  In this experiment, we 'oversample' the input signal by a
  factor 2 so that $m = 2n$.  (a) Success rate of $(P_1)$ for $m =
  512$.  (b) Success rate of $(P_1)$ for $m = 1024$. Observe the
  similar pattern and cut-off point. In these experiments, exact
  recovery occurs as long as about 17\% or less of the entries are
  corrupted.}
\label{fig:reccurves2}
\end{figure}

\begin{figure}
\centerline{
\begin{tabular}{c}
\includegraphics[width=3in]{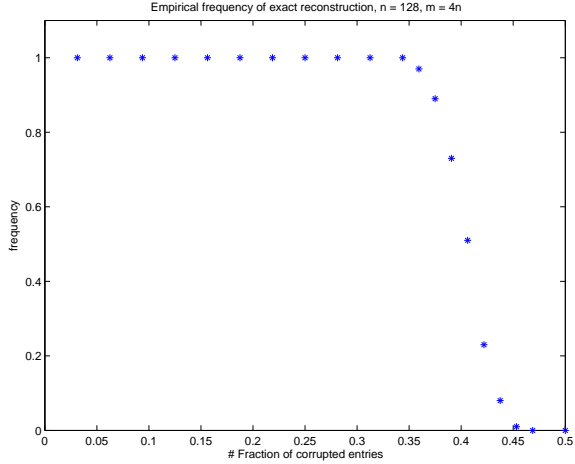} 
\end{tabular}
}
\caption{\small\sl  $\ell_1$-recovery of an input signal from 
  $y = Af + e$ with $A$ an $m$ by $n$ matrix with independent Gaussian
  entries.  In this experiment, we 'oversample' the input signal by a
  factor 4 so that $m = 4n$.  In these experiments, exact recovery
  occurs as long as about 34\% or less of the entries are corrupted.}
\label{fig:reccurves4}
\end{figure}

The results are presented in Figure~\ref{fig:reccurves2} and
Figure~\ref{fig:reccurves4}.  Figure~\ref{fig:reccurves2} examines the
situation in which the length of the code is twice that of the input
vector $m = 2n$, for $m = 512$ and $m = 1024$. Our experiments show
that one recovers the input vector {\em all the time} as long as the
fraction of the corrupted entries is below 17\%. This holds for $m =
512$ (Figure~\ref{fig:reccurves2}(a)) and $m = 1024$
(Figure~\ref{fig:reccurves2}(b)).  In Figure~\ref{fig:reccurves4}, we
investigate how these results change as the length of the codewords 
increases compared to the length of the input, and examine the
situation in which $m = 4n$, with $m = 512$.  Our experiments show
that one recovers the input vector {\em all the time} as long as the
fraction of the corrupted entries is below 34\%.

\section{Optimal Signal Recovery}
\label{sec:optimalrecovery}

Our recent work \cite{OptimalRecovery} developed a set of ideas
showing that it is surprisingly possible to reconstruct interesting
classes of signals accurately from highly incomplete measurements. The
results in this paper are inspired and improve upon this earlier work
and we now elaborate on this connection. Suppose we wish to
reconstruct an object $\alpha$ in $\R^m$ from the $K$ linear
measurements
\begin{equation}
\label{eq:linmeas}
y_k = \<\alpha,\phi_k\> \quad k = 1,\ldots,K \qquad
\text{or}\quad\quad y = F \alpha, 
\end{equation}
with $\phi_k$, the $k$th row of the matrix $F$.  Of special interest
is the vastly underdetermined case, $K << N$, where there are many
more unknowns than observations. We choose to formulate the problem
abstractly but for concreteness, we might think of $\alpha$ as the
coefficients $\alpha = \Psi^* f$ of a digital signal or image $f$ in
some nice orthobasis, e.g. a wavelet basis so that the information
about the signal is of the form $y = F \alpha = F \Psi^* f$.

Suppose now that the object of interest is {\em compressible} in the
sense that the reordered entries of $\alpha$ decay like a power-law;
concretely, suppose that the entries of $\alpha$, rearranged in
decreasing order of magnitude, $|\alpha|_{(1)} \geq |\alpha|_{(2)} \geq
\cdots \geq |\alpha|_{(m)}$, obey 
\begin{equation}
\label{eq:wlp}
|\alpha|_{(k)} \leq B\cdot k^{-s}
\end{equation}
for some $s \geq 1$.  We will denote by $\cF_s(B)$ the class of all
signals $\alpha \in \R^m$ obeying \eqref{eq:wlp}. The claim is that it is
possible to reconstruct compressible signals from only a small number
of random measurements.
\begin{theorem}
\label{teo:optimalrecovery}
Let $F$ be the measurement matrix as in \eqref{eq:linmeas} and consider
the solution $\alpha^\sharp$ to
  \begin{equation}
    \label{eq:newP1}
\min_{\tilde \alpha \in \R^m} \, \|\tilde \alpha\|_{\ell_1} 
\quad \text{ subject to } \quad F \tilde \alpha = y.
  \end{equation}
  Let $S \le K$ such that $\delta_S + 2\theta_S + \theta_{S,2S} < 1$
  and set $\lambda = K/S$.  Then $\alpha^\sharp$ obeys
  \begin{equation}
    \label{eq:optimalrecovery}
   \sup_{\alpha \in \cF_s(B)}  \|\alpha - \alpha^\sharp\| \le C 
\cdot (K/\lambda)^{-(s-1/2)}.  
  \end{equation}
\end{theorem}
To appreciate the content of the theorem, suppose one would have
available an {\em oracle} letting us know which coefficients
$\alpha_k$, $1 \le k \le m$, are large (e.g. in the scenario we
considered earlier, the oracle would tell us which wavelet
coefficients of $f$ are large). Then we would acquire information
about the $K$ largest coefficients and obtain a truncated version
$\alpha_K$ of $\alpha$ obeying 
\[
\|\alpha - \alpha_K\| \ge c \cdot K^{-(s-1/2)},  
\] 
for generic elements taken from $\cF_s(B)$. Now
\eqref{eq:optimalrecovery} says that not knowing anything about the
location of the largest coefficients, one can essentially obtain the
same approximation error by nonadaptive sampling, provided the number
of measurements be increased by a factor $\lambda$. The larger $S$,
the smaller the oversampling factor, and hence the connection with the
decoding problem. Such considerations make clear that Theorem
\ref{teo:optimalrecovery} supplies a very concrete methodology for
recovering a compressible object from limited measurements and as
such, it may have a significant bearing on many fields of science and
technology. We refer the reader to \cite{OptimalRecovery} and
\cite{CompressedSensing} for a discussion of its implications.

Suppose for example that $F$ is a Gaussian random matrix as in Section
\ref{sec:random}. We will assume the same special normalization so
that the variance of each individual entry is equal to $1/K$.
Calculations identical to those from Section \ref{sec:random} give
that with overwhelming probability, $F$ obeys the hypothesis of the
theorem provided that
\[
S \le K/\lambda, \quad \lambda = \rho \cdot \log(m/K), 
\]
for some positive constant $\rho > 0$. Now consider the statement of
the theorem; there is a way to invoke linear programming and obtain a
reconstruction based upon $O(K \log (m/K))$ measurements only, which
is at least as good as that one would achieve by knowing all the
information about $f$ and selecting the $K$ largest coefficients. In
fact, this is an optimal statement as \eqref{eq:optimalrecovery}
correctly identifies the minimum number of measurements needed to
obtain a given precision. In short, it is impossible to obtain a
precision of about $K^{-(s-1/2)}$ with fewer than $K \log (m/K)$
measurements, see \cite{OptimalRecovery,CompressedSensing}.

Theorem \ref{teo:optimalrecovery} is stronger than our former result,
namely, Theorem 1.4 in \cite{OptimalRecovery}. To see why this is
true, recall the former claim: \cite{OptimalRecovery} introduced two
conditions, the uniform uncertainty principle (UUP) and the exact
reconstruction principle (ERP).  In a nutshell, a {\em random} matrix
$F$ obeys the UUP with oversampling factor $\lambda$ if $F$ obeys
\begin{equation}
  \label{eq:uup}
   \delta_S \le 1/2, \qquad S = \rho \cdot K/\lambda, 
\end{equation}
with probability at least $1 - O(N^{-\gamma/\rho})$ for some fixed
positive constant $\gamma > 0$. Second, a measurement matrix $F$ obeys
the ERP with oversampling factor $\lambda$ if for each fixed subset
$T$ of size $|T| \le S$ \eqref{eq:uup} and each `sign' vector $c$
defined on $T$, there exists with the same overwhelmingly large
probability a vector $w \in H$ with the following properties:
\begin{itemize}
\item[(i)] $\<w, v_j\> = c_j$, for all $j \in T$;
\item[(ii)] and $|\<w, v_j\>| \leq \frac{1}{2}$ 
for all $j$ not  in $T$. 
\end{itemize}
Note that these are the conditions listed at the beginning of section
\ref{sec:main} except for the 1/2 factor on the complement of $T$. Fix
$\alpha \in \cF_s(B)$.  \cite{OptimalRecovery} argued that if a
random matrix obeyed the UUP and the ERP both with oversampling factor
$\lambda$, then 
\[
\|\alpha - \alpha^\sharp\| \le C \cdot (K/\lambda)^{-(s-1/2)},
\]
with inequality holding with the same probability as before.  Against
this background, several comments are now in order:
\begin{itemize}
\item First, the new statement is more general as it applies to all
  matrices, not just random matrices. 
  
\item Second, whereas our previous statement argued that for each
  $\alpha \in \R^m$, one would be able---with high probability---to
  reconstruct $\alpha$ accurately, it did not say anything about the
  worst case error for a fixed measurement matrix $F$. This is an
  instance where the order of the quantifiers plays a role. Do we need
  different $F$'s for different objects? Theorem
  \ref{teo:optimalrecovery} answers this question unambiguously; the
  same $F$ will provide an optimal reconstruction for {\em all} the
  objects in the class.
  
\item Third, Theorem \ref{teo:optimalrecovery} says that the ERP
  condition is redundant, and hence the hypothesis may be easier to
  check in practice.  In addition, eliminating the ERP isolates the
  real reason for success as it ties everything down to the UUP. In
  short, the ability to recover an object from limited measurements
  depends on how close $F$ is to an orthonormal system, but only when
  restricting attention to sparse linear combinations of columns.
\end{itemize}

We will not prove this theorem as this is a minor modification of that
of Theorem 1.4 in the aforementioned reference. The key point is to
observe that if $F$ obeys the hypothesis of our theorem, then by
definition $F$ obeys the UUP with probability one, but $F$ also obeys
the ERP, again with probability one, as this is the content of Lemma
\ref{erp-3}.  Hence both the UUP and ERP hold and therefore, the
conclusion of Theorem \ref{teo:optimalrecovery} follows. (The fact that
the ERP actually holds for all sign vectors of size less than $S$
is the reason why \eqref{eq:optimalrecovery} holds uniformly over all
elements taken from $\cF_s(B)$, see \cite{OptimalRecovery}.)

\section{Discussion}
\label{sec:discussion}

\subsection{Connections with other works}

In our linear programming model, the plaintext and ciphertext had
real-valued components.  Another intensively studied model occurs when
the plaintext and ciphertext take values in the finite field $F_2 :=
\{0,1\}$.  In recent work of Feldman et al. \cite{feldman-thesis},
\cite{feldman-1}, \cite{feldman-2}, linear programming methods (based
on relaxing the space of codewords to a convex polytope) were
developed to establish a polynomial-time decoder which can correct a
constant fraction of errors, and also achieve the
information-theoretic capacity of the code.  There is thus some
intriguing parallels between those works and the results in this
paper, however there appears to be no direct overlap as our methods
are restricted to real-valued texts, and the work cited above requires
texts in $F_2$.  Also, our error analysis is deterministic and is thus
guaranteed to correct arbitrary errors provided that they are
sufficiently sparse.

The ideas presented in this paper may be adapted to recover input
vectors taking values from a finite alphabet. We hope to report on
work in progress in a follow-up paper.

\subsection{Improvements}

There is little doubt that more elaborate arguments will yield
versions of Theorem \ref{teo:gauss} with tighter bounds.  Immediately
following the proof of Lemma \ref{erp-3}, we already remarked that one
might slightly improve the condition $\delta_S + \theta_{S,S} +
\theta_{S,2S} < 1$ at the expense of considerable complications. More
to the point, we must admit that we used well-established tools from
Random Matrix Theory and it is likely that more sophisticated ideas
might be deployed successfully. We now discuss some of these.

Our main hypothesis reads $\delta_S + \theta_{S,S} + \theta_{S,2S} <
1$ but in order to reduce the problem to the study of those
$\delta$ numbers (and use known results), our analysis actually relied
upon the more stringent condition $\delta_S + \delta_{2S} +
\delta_{3S} < 1$ instead, since
\[
\delta_S + \theta_{S,S} + \theta_{S,2S} \le \delta_S + \delta_{2S} +
\delta_{3S}. 
\]
This introduces a gap. Consider a fixed set $T$ of size
$|T| = S$. Using the notations of
that Section \ref{sec:random}, we argued that 
\[
\delta(F_T) \approx 2\sqrt{S/p} + S/p, 
\] 
and developed a large deviation bound to quantify the departure from
the right hand-side. Now let $T$ and $T'$ be two disjoint sets of
respective sizes $S$ and $S'$ and consider $\theta(F_T, F_{T'})$:
$\theta(F_T, F_{T'})$ is the cosine of the principal angle between the
two random subspaces spanned by the columns of $F_T$ and $F_{T'}$
respectively; formally
\[ 
\theta(F_T, F_{T'}) = \sup \<u, u'\>, 
\qquad u \in \text{span}(F_T), \, u' \in \text{span}(F_{T'}), \, 
\|u\| = \|u'\| = 1.
\] 
We remark that this quantity plays an important analysis in
statistical analysis because of its use to test the significance of
correlations between two sets of measurements, compare the literature
on {\em Canonical Correlation Analysis} \cite{Muirhead}. Among other
things, it is known \cite{Wachter} that 
\[
\theta(F_T, F_{T'}) \goto \sqrt{\gamma(1-\gamma')} + \sqrt{\gamma'(1-\gamma)} 
\quad \text{a.s.} 
\] 
as $p \goto \infty$ with $S/p \goto \gamma$ and $S'/p \goto \gamma'$.
In other words, whereas we used the limiting behaviors
\[
\delta(F_{2T}) \goto 2\sqrt{2\gamma} + 2\gamma, \quad \delta(F_{3T}) \goto  
2\sqrt{3\gamma} + 3\gamma,
\]
there is a chance one might employ instead 
\[ 
\theta(F_T, F_{T'}) \goto 2\sqrt{\gamma(1-\gamma)}, \quad \theta(F_T, F_{T'})
\goto \sqrt{\gamma(1-2\gamma)} + \sqrt{2\gamma(1-\gamma)}
\] 
for $|T| = |T'| = S$ and $|T'| = 2|T| = 2S$ respectively, which is
better. Just as in Section \ref{sec:random}, one might then look for
concentration inequalities transforming this limiting behavior into
corresponding large deviation inequalities.  We are aware of very
recent work of Johnstone and his colleagues \cite{IainWaldWeb} which
might be here of substantial help.

Finally, tighter large deviation bounds might exist together with more
clever strategies to derive uniform bounds (valid for all $T$ of size
less than $S$) from individual bounds (valid for a single $T$).  With
this in mind, it is interesting to note that our approach hits a limit
as
\begin{equation}
  \label{eq:lower-limit}
  \liminf_{S \goto \infty, \, S/m \goto r} \,\, 
\delta_S + \theta_{S,S} + \theta_{S,2S}
\ge J(m/p \cdot r),
\end{equation}
where $J(r) := 2\sqrt{r} + r + (2 + \sqrt{2})\sqrt{r(1-r)} +
\sqrt{r(1-2r)}$.  Since $J(r)$ is greater than 1 if and only if $r >
2.36$, one would certainly need new ideas to improve Theorem
\ref{teo:gauss} beyond cut-off point in the range of about 2\%. The
lower limit \eqref{eq:lower-limit} is probably not sharp since it does
not explicitly take into account the ratio between $m$ and $p$; at
best, it might serve as an indication of the limiting behavior when
the ration $p/m$ is not too small.

\subsection{Other coding matrices}

This paper introduced general results stating that it is possible to
correct for errors by $\ell_1$-minimization. We then explained how the
results specialize in the case where the coding matrix $A$ is sampled
from the Gaussian ensemble. It is clear, however, that one could use
other matrices and still obtain similar results; namely, that $(P'_1)$
recovers $f$ exactly provided that the number of corrupted entries
does not exceed $\rho \cdot m$.  In fact, our previous work suggests
that partial Fourier matrices would enjoy similar properties
\cite{RandomBP,OptimalRecovery}. Other candidates might be the
so-called {\em noiselets} of Coifman, Geshwind and Meyer
\cite{Noiselets}. These alternative might be of great practical
interest because they would come with fast algorithms for applying $A$
or $A^*$ to an arbitrary vector $g$ and, hence, speed up the
computations to find the $\ell_1$-minimizer.

\end{document}